\documentclass{amsart}
\usepackage{subfigure}
\usepackage{graphicx}
\usepackage{psfrag}

\usepackage{amsmath}
\usepackage{amssymb}
\usepackage{caption}

\newcommand{\1}{\mathbf 1}
\newcommand{\dd}{\mathrm d}
\newcommand{\e}{\mathrm e}
\newcommand{\topo}{{\mbox{\tiny \rm top}}}

\renewcommand{\phi}{\varphi}

\renewcommand{\tilde}{\widetilde}

\newtheorem{theorem}{Theorem}[section]
\newtheorem{lemma}[theorem]{Lemma}
\newtheorem{proposition}[theorem]{Proposition}
\newtheorem{corollary}[theorem]{Corollary}
\newtheorem{definition}[theorem]{Definition}

\def\R{\mathbb{R}}
\def\C{\mathbb{C}}
\def\N{\mathbb{N}}

\def\Sp{\mathbb{S}}
\def\S{\Sigma_{d}}

\def\a{\alpha}

\def\d{\delta}
\def\c{\ell}
\def\chi{\ell}

\begin{document}
\author{Marc Kesseb\"ohmer and Bernd O. Stratmann}
\address{Fachbereich 3 - Mathematik und Informatik, Universit\"at Bremen, D--28359 Bremen, Germany}
	 \email{mhk@math.uni-bremen.de}	 
\address{Mathematical Institute, University of St Andrews, St Andrews KY16 9SS, Scotland}
         \email{bos@maths.st-and.ac.uk}

\title[Refined measurable rigidity and flexibility]{
Refined measurable rigidity and flexibility for conformal
iterated function systems}


\date{\today}


\begin{abstract}
In this paper we investigate aspects of rigidity and flexibility 
for conformal
iterated function systems. For the case in which the systems are not essentially affine
we show that two such systems are conformal equivalent if and only 
if in each of their Lyapunov spectra there exists at least one level 
set such that the corresponding Gibbs measures coincide. We then 
proceed by 
comparing this result with the essentially affine situation. We show 
that   essentially affine systems are far less rigid than non--essentially
affine systems, and  
subsequently we then investigate  the extent  of their flexibility.
\end{abstract}
 \maketitle

\section{Introduction}

In 1982 D. Sullivan published his influential purely measurable
form of Mostow's rigidity theorem. It states that if two
geometrically finite Kleinian groups are conjugate under a Borel map 
$F$
which is non-singular with respect to  the  Patterson measures 
associated with the two groups, then $F$ agrees almost everywhere
with a conformal conjugacy (\cite{sul1}, see also \cite{sul2},
\cite{sul3} and \cite{B}). Since the appearance of this 
theorem  the concept of measurable 
rigidity has attracted a great deal of attention, and in the 
meanwhile numerous 
generalisations and variations have been obtained. One of these 
 was derived by Hanus and Urba\'nski (\cite{HU}),
 who considered   non-essentially affine, conformal
iterated function systems (see Section \ref{section2} for the 
definitions),
and showed that two such systems $\Phi$ and $\Psi$ are
conformal equivalent
if and only if their associated conformal measures $\mu_{\Phi}$ and 
$\mu_{\Psi}$ (each of  
maximal Hausdorff dimension) coincide up to permutation of the 
generators.  This result can be seen  as
 the 
 starting point for this paper.

Our first goal is to  give a multifractal refinement 
of the result in \cite{HU}, where for ease of exposition we restrict the discussion to the
$1$--dimensional finite case. 
For this we will recall that  each system  $\Phi$ 
gives rise to its Lyapunov spectrum $u\mapsto \c_{\Phi}(u)$,  
which is given by the 
multifractal spectrum of the  measure of maximal entropy associated 
with $\Phi$. 
Moreover, each  level set in this spectrum supports 
a canonical 
shift--invariant Gibbs measure $\mu_{\Phi, u}$.  
In a nut shell, our main result for non-essentially affine, conformal
iterated function systems  is that two such systems 
$\Phi$ and $\Psi$ are 
conformal equivalent if and only if $\mu_{\Phi,u}$ is equal to $\mu_{\Psi,v}$ 
up to permutation of the generators, for some
$u,v \in \R\setminus \{0\}$
(see Theorem \ref{NAS} for a more complete statement which also 
involves 
cohomological equivalence of the 
associated canonical geometric potential functions, equality of pressure functions as 
well as equality of Lyapunov spectra).
 
In the second part of the paper we consider essentially 
affine, 
 conformal iterated function systems. 
  Note that for non-essentially affine systems  a conjugation 
 map 
 between two systems is conformal if and only if it is bi-Lipschitz 
 (see \cite{MU} Theorem 7.2.4). Hence, for essentially affine
 systems bi-Lipschitz conjugation is the natural substitute for 
 conformal conjugation. 
 By 
investigating similar questions as before for the non-essentially 
affine case, we obtain that 
from the point of view of multifractal rigidity essentially 
affine systems  
behave  rather different  than non-essentially 
affine systems. 
For instance,  if  for two essentially 
affine systems $\Phi$ and $\Psi$  we have that $\mu_{\Phi,u}$ is equal to $\mu_{\Psi,v}$ 
up to permutation of the generators, for 
some $u,v \in \R\setminus \{0\}$, then 
this does {\em not necessarily}  imply that $\Phi$ and $\Psi$ are 
bi-Lipschitz equivalent. More precisely, we show that 
equality of $\mu_{\Phi,u}$ and $\mu_{\Psi,u}$ up to permutation 
of the generators {\em together} with the equality of the pressure functions 
$P_{\Phi}$ and $P_{\Psi}$ at $u$, for some $u \in \R \setminus \{0\}$,
is equivalent to the fact that $\Phi$ is bi-Lipschitz equivalent 
to $\Psi$, as well as to the facts $P_{\Phi}= P_{\Psi}$, 
$\c_{\Phi}=\c_{\Psi}$ and cohomological equivalence of the two canonical geometric potential functions 
associated with the systems (see Theorem \ref{linearrigid}).  
These results clearly show that essentially affine systems are less
rigid than non-essentially affine systems, and a further 
investigation of this phenomenon of flexibility  is then given
in Section \ref{Subsec:MultiFlex}.
There, we derive sufficient and necessary conditions for equality
of $\mu_{\Phi,u}$ and $\mu_{\Psi,v}$ in terms of the pressure 
functions and the canonical geometric potential functions 
(see Theorem \ref{Thm:Multiflex}).
Also, we show that this situation does in fact occur. Namely,
in Proposition \ref{flexsection} we obtain that if $\mu_{\Phi,u}$ 
is given and $v$ fulfils a 
certain admissibility condition (see Definition \ref{adm}), then
there exists an essentially 
affine system $\Psi$ such that $\mu_{\Phi,u}$ is equal to 
$\mu_{\Psi,v}$ up to permutation of the generators. Finally, we  
 give a 
brief discussion of  the extent of flexibility  of an 
essentially affine system. The outcome here is  
 that for a
non-degenerate $\Phi$ the set of systems $\Psi$ for which 
$\mu_{\Phi,u}$ is equal to $\mu_{\Psi,v}$ up to permutation of the 
generators, for some $u,v\in
\mathbb{R} \setminus \{0\}$,
forms a $2$--dimensional submanifold of 
the moduli space of $\Phi$, whereas if $\Phi$ is 
degenerate 
then this set is a $1$--dimensional submanifold (see Proposition 
\ref{ending}).

\section{Preliminaries}\label{section2}
\subsection{Conformal iterated function systems}
Throughout this paper we consider  \emph{conformal iterated function
systems} (CS) on
some connected
compact set $X \subset \R$.
Recall from \cite{HU}  (see also \cite{MU}) that these systems are  generated by 
an ordered family 
$\Phi$ of injective contractions $(\phi_i : X \to \hbox{Int} X \, |  \, i \in 
I)$, for
some given finite index set $I:=\{1,\ldots,d\}$ with at 
least two elements. Furthermore, $\Phi$  satisfies the 
following conditions, where we use the notation $\phi_\omega :=  \phi_{x_1}\circ \phi_{x_2}\circ \ldots 
\circ \phi_{x_n} $ for  $\omega = x_1 x_2 \ldots x_n \in I^n$.
\begin{description}
\item[{\it Open set condition}]  $\phi_{i} (\hbox{Int} (X) ) \subset 
\hbox{Int} (X)$ for all $i \in I$, \\ and 
$\phi_i (\hbox{Int} (X)) 
\cap \phi_j (\hbox{Int}(X))= \emptyset  \,\hbox{ for each }  \, i, j \in I, 
i\neq j$.
\item[{\it Conformality--condition}] There exists an open 
connected set $U \subset \C$ containing $X$ such that $\phi_{i}$
extends to a conformal map on $U$, for each $i \in I$.
\item[{\it Bounded distortion property}] There exists $C \ge 1$ 
such that for all $n \in \N,
 \omega  \in I^n$ and $x,y \in U$ we have 
$$
|\phi_\omega'(y)|  \le C \;  |\phi_\omega'(x)|.
$$
\end{description}
A central object associated with a CS $\Phi$  is its {\em limit set} 
  \[\Lambda(\Phi):=\bigcap_{n\in\N}\bigcup_{\omega\in I^n} \phi_\omega (X).\]
Clearly, $\Lambda (\Phi)$ is the unique non-empty compact subset 
of $\R$ for which
$\Lambda(\Phi) = \bigcup_{i\in I} \phi_{i}(\Lambda(\Phi))$.
From a combinatorial point of view  $\Phi$ is
described by the 
full-shift $\Sigma_{d} := I^{{\N}}$.
As usual, we assume $\Sigma_{d}$ to be equipped with the left-shift 
map $\sigma$.  The link between $\S$ and $\Phi $
is provided by the canonical bijection $\pi_{\Phi}: \S \to \Lambda(\Phi)$ 
which is given by $\pi_{\Phi}(x_{1} x_{2} \ldots) := \lim_{n \to 
\infty} \phi_{x_{1} x_{2}  \ldots x_{n}} (X)$.   
Evidently,  we can always think of  $\Phi $ as being a  
conformal fractal
representation
  of $\Sigma_{d}$. \\
 We also  
 consider the special situation in which all the $\phi_{i}$
 are in particular affine 
   transformations. In this case the system is called
    {\em affine iterated function system} (AS), and 
occasionally it will also be referred to as an {\em 
affine fractal representation} of 
  $\Sigma_{d}$.   One of the major issues of this paper is to study 
  certain deformations of a given CS  $\Phi= (\phi_i : X \to \hbox{Int} X \, |  \, i \in 
I)$. 
  More precisely,   
  let     
 $\Psi:=(\psi_i:Y\to \hbox{Int}Y \, | \, i\in I)$ be 
 some other CS defined on some connected compact set 
 $Y \subset \R$. 
 Then $\Psi$  is called a {\em deformation} of 
  $\Phi $  if 
there  exists 
a  bi-Lipschitz map $h:\Lambda(\Phi)\to \Lambda(\Psi)$
such that 
\[ \psi_i= h\circ \phi_i \circ h^{-1},\hbox{  for each } i\in I.\] 
A map $h$ of this type
will be called 
 a {\em fractal boundary 
correspondence}.
In particular,  if in here $\Phi$ is an AS, that is if $\Psi$ is a deformation 
of an affine iterated function system, then $\Psi$ will
be referred to as {\em essentially affine iterated function 
system} (EAS). On the other hand, if $\Phi$ is a CS which is not an EAS
then $\Phi$  will be called {\em non-essentially affine iterated function 
system} (NAS).
\\ Let us also introduce the
 {\em deformation space} 
${\mathcal T}(\Sigma_{d})$ associated with  
$\Sigma_{d}$. This is given by
 \[{\mathcal T}(\Sigma_{d}):= \{ \Psi: \Psi   \hbox{ is a 
 CS on $\Sigma_{d}$} \}
.\]
Clearly,  ${\mathcal T}(\Sigma_{d})$ relates to $\Sigma_{d}$ 
similar as 
the Teichm\"uller space for a Riemann surface relates to the 
associated fundamental group. 
We then decompose  the space  ${\mathcal T}(\Sigma_{d}) $
    into the two disjoint deformation spaces
    \[{\mathcal T}_{E}(\Sigma_{d}):= \{ \Psi: \Psi   \hbox{ is an 
     EAS on $\Sigma_{d}$} \} \hbox{  and  } {\mathcal T}_{N}(\Sigma_{d}):= \{ \Psi: \Psi   \hbox{ is an 
NAS on $\Sigma_{d}$} \}.\]
Also, we introduce an equivalence relation 
 on ${\mathcal T} 
(\Sigma_{d})$ as follows. Two systems 
$\Phi, \Psi \in {\mathcal T}(\Sigma_{d})$ are said to be
equivalent  ($\Phi  
\sim \Psi$)   if and only 
if there exists
   a fractal boundary correspondence $h: \Lambda(\Phi) \to 
   \Lambda(\Psi)$ between them. 
  Finally, recall that a CS is called {\em degenerate} if it is 
      equivalent to an AS $\Psi= (\psi_i : X \to \hbox{Int} X \, |  \, i \in 
I)$ for which $\psi_{i}' = \psi_{j}'$, for all $i,j \in I$. 
 It is easy to see 
 that for a degenerate EAS  the multifractal analysis in this paper is trivial.

  \subsection{Thermodynamic and multifractal formalism for CS}
  Let
  $\Phi =(\phi_i : X \to \hbox{Int} X \, |  \, i \in 
      I) \in {\mathcal T}(\Sigma_{d}) $ be  given, and let $\delta_{\Phi}$ refer to the
      Hausdorff dimension of  $\Lambda(\Phi)$. 
  Throughout, we require the following  standard concepts  
  from thermodynamic 
  formalism, and we assume that the reader is familiar with the basics of 
  this formalism (see e.g. \cite{Bowen}, \cite{Denker}, \cite{Pesin}, 
  \cite{Ruelle}). Here we use the common notation $[x_{1}\ldots x_{n}]:=
  \{y=(y_1 y_2 \ldots)  \in \Sigma_{d}: y_{i}=x_{i} \hbox{ for } 
  i=1,\ldots,n\}$ and $S_{n} f:= \sum_{k=0}^{n-1} f \circ \sigma^{k}$.
     \begin{itemize}
      \item  The  {\em canonical geometric potential} $I_{\Phi}:\S \to \R$
      associated with $\Phi$  is given 
      by   $I_{\Phi}(x):=\log 
  \phi_{x_1}'(\pi_{\Phi}(x))$ for all $x=(x_1 x_2 \ldots) 
  \in\Sigma_d$.
  \item   $\mu_{\Phi}$ refers to a 
  Gibbs measure on $\S$ for the potential 
  $\d_{\Phi} I_{\Phi}$.
  \item  $\c_{\Phi}$ refers to the {\em Lyapunov spectrum}  of $\Phi$, 
   given for $\alpha 
   \in \R$ by 
   \[ \c_{\Phi}(\alpha):=\dim_{H}\left(\pi_{\Phi} \left(\left\{x \in 
   \Sigma_{d}: 
     \lim_{n \to 
   \infty} \frac{S_n  I_\Phi(x)}{-n} = \alpha  
   \right\}\right) \right)  . \]
  \item 
  $\mathcal{P}_{\Phi} : \R \to \R$  denotes the {\em pressure 
  function},
  given
  for a potential function  $f:\Sigma_{d} \to \R$ by 
  \[\mathcal{P}_\Phi (f):=\lim _{n\to \infty}\frac{1}{n}\log\sum_{\omega\in 
  I^n}\exp( \sup_{x\in [\omega]} S_n  f (x)) .
  \]
  Also,    for $u 
  \in \R$ we define
  \[ P_{\Phi} (u) : = \mathcal{P}_{\Phi} (u \, I_{\Phi}) \]
  and \[
   \alpha_{\Phi}(u) := - P_{\Phi}'(u) =- \int I_{\Phi} \, \dd 
  \mu_{\Phi,u} .\]
  Here, $\mu_{\Phi,u}$ refers to the {\em $\sigma$--invariant 
    Gibbs measure} on $\S$ for the potential function
    $u I_{\Phi} -P(u)$, where `Gibbs' means as usual  that 
    for all  $n \in \N, (x_{1},\ldots,x_{n}) \in I^{n}$ and 
    $x \in [x_{1} \ldots x_{n}]$,
    \[ \mu_{\Phi,u} ([x_{1},\ldots,x_{n}]) \asymp e^{u  S_{n} 
    I_{\Phi}(x)  - n P_{\Phi}(u)} .\]
    Furthermore, we let $m_{\Phi,u}$ denote the {\em $(u I_{\Phi} 
    -P_{\Phi}(u))$--conformal measure} within the measure class of 
    $\mu_{\Phi,u}$, given by
    \[ \frac{d m_{\Phi,u} \circ \phi_{i}}{d 
	       m_{\Phi,u}} = \left| \phi'_{i}\right|^{u}  \, 
	       \e^{-P_{\Phi}(u)} , \hbox{ for all } i \in I.
    \]
  \end{itemize}
  Finally, throughout 
   we require the following notions of equivalence in 
  connection with
   two given $\Phi, \Psi \in \mathcal{T}(\Sigma_{d})$.  
  \begin{itemize}
    \item  The Gibbs measures $\mu_{\Phi,u}$ 
  and $\mu_{\Psi,v}$ are {\em equal up 
  to permutation} ($\mu_{\Phi,u} \cong \mu_{\Psi,v}$) if and only if 
  $\mu_{\Phi,u} = \mu_{\Psi_{0},v}$ for some system $\Psi_{0}$ obtained from $\Psi$
  by a permutation of the generators of $\Psi$. 
  \item The potentials $I_{\Phi}$ 
  and $I_{\Psi}$ are {\em cohomological equivalent} ($I_{\Phi} \simeq I_{\Psi}$) if 
  $I_{\Phi}$ is cohomologous to $I_{\Psi_{0}}$, for some system $\Psi_{0}$ 
  obtained from $\Psi$
  by a permutation of the generators of $\Psi$. (Recall that two 
  functions
  $f,g: \S \to \R$ are {\em cohomologous} if there exists a continuous 
  function $e:\S \to \R$ such that $f-g=e- e\circ \sigma$). 
 \end{itemize}
  
For the type of  iterated function systems which we consider in this
paper the calculation of the Lyapunov spectrum is basically  an
application of the multifractal analysis of the 
measure of maximal entropy for cookie-cutter Cantor sets. The 
following proposition summarises the outcome of this analysis. 
Subsequently, we 
will outline the proof employing the down-to-earth approach
given in \cite{Falconer}. Note that the proposition also  
immediately  follows from the multifractal formalism for growth rates
developed in \cite{KS1}. Also, note that in here the function 
$\beta_{\Phi}$ is precisely the inverse of the function $\alpha_{\Phi}$,
that is $\beta_{\Phi} \circ \alpha_{\Phi}= id.$.

\begin{proposition}\label{MFF}
      Let $\Phi \in {\mathcal T}(\Sigma_{d}) $ non-degenerate 
      be given. Then there exists
      a real-analytic function $\beta_{\Phi}: \R \to \R$ and
      $\alpha_{-}, \alpha_{+} >0$ such that $\c_{\Phi}(\alpha)
      =0$ for all $ \alpha \notin 
      (\alpha_{-},\alpha_{+})$ and such that for all $ \alpha \in 
      [\alpha_{-},\alpha_{+}]$,
      \[ \c_{\Phi}(\alpha)  =  
      \beta_{\Phi}(\alpha) 
      +\frac{P_{\Phi}\left(\beta_{\Phi}(\alpha)\right)}{\alpha}.\]
\end{proposition}
 \begin{proof} (Sketch) Let $\nu$ refer to the measure of maximal 
 entropy for the system $\Phi$ on $\Sigma_{d}$. Then $\nu$ is a Gibbs measure
 for the potential function $\varphi$ constant equal to the negative 
 of the topological entropy 
 $h_\topo:=\log d$.
 Hence, we in particular have 
 $\nu([\omega]) \asymp \exp(S_{n} \varphi 
 (x)) = d^{-n}$ for all $n \in \N, \omega 
 \in I^{n}$ and $x \in [\omega]$. Trivially, we have $\varphi < 0$ and
 $\mathcal{P}(\varphi)=0$, which shows that $\nu$ can be 
 analysed by standard multifractal analysis (see e.g.  \cite{Falconer}). This 
 gives that there exists a well-defined, strictly decreasing, 
 real-analytic function $\gamma_{\Phi}:\R \to \R$ such that 
 $\mathcal{P}(\gamma_{\Phi}(t) \, I_{\Phi}+ t \varphi) =0$, for all 
 $t \in \R$. In order to determine the Hausdorff dimension spectrum
 of \[ E_{\tau}:= \left\{x \in 
 \Lambda(\Phi) : \lim_{n\to \infty} \frac{\log 
 \nu\left(\pi_{\Phi}^{-1} (B(x,r))\right)}{\log r}  = \tau
 \right\}  ,\]
 one considers the Legendre transform of $\gamma_{\Phi}$, given by 
 $f(\tau) =\inf \{ \gamma_{\Phi} (t) + t \, \tau: t \in \R\}$, or what is 
 equivalent $f(\tau) = \gamma_{\Phi} (t_{\tau}) + t_{\tau} \, \tau$ where $t_{\tau}$
 is determined by $\gamma_{\Phi}'(t_{\tau}) =- \tau$. In particular, 
 there exists a maximal interval $(\tau_{-}, \tau_{+})$ on which $f$ 
 is continuous, concave and strictly positive; outside  this 
 interval $f$ vanishes. Now, the key 
 observation is that  there exists a Gibbs measure 
 $\nu_{\tau}$
 for the potential function $\gamma_{\Phi}(t_{\tau})\, I_{\Phi}+ 
 t_{\tau} \varphi$ which is concentrated on $\pi_{\Phi}^{-1}(E_{\tau})$.
 (Note that the measure $\nu_{\tau}$ coincides with the measure 
 $\mu_{\Phi,\gamma_{\Phi}(t_{\tau})}$ which we  already introduced 
 above).
 Hence, we have for all $n \in \N$, $\omega \in I^{n}$ and $x 
 \in [\omega]$,
 \[ \nu_{\tau} ([\omega])  \asymp \exp 
 \left(\gamma_{\Phi}(t_{\tau})\, S_{n}I_{\Phi}(x)+ 
 t_{\tau}  S_{n}\varphi (x)\right) .\]
 Since by the bounded distortion property $\exp(S_{n}I_{\Phi}(x)) \asymp 
|\pi_{\Phi}([\omega]) |$, the mass distribution 
 principle   therefore immediately gives $\dim_{H}(E_{\tau}) = f(\tau)$. 
 To 
 finish the proof, note that 
 \begin{eqnarray*}  \pi_{\Phi}^{-1}(E_{\tau}) &=& \left\{ x=(x_{1}x_{2} \ldots)\in 
 \Sigma_{d} : \lim_{n \to \infty} \frac{-n \, h_{\topo}}{\log 
 |\pi_{\Phi}([x_{1}\ldots x_{n}])|} 
  = \tau\right\} \\  &=& 
 \left\{ x \in 
  \Sigma_{d} : \lim_{n \to \infty} \frac{\log S_{n} I_{\Phi}(x)}{-n} = 
  \frac{h_\topo}{\tau} \right\} .\end{eqnarray*}
  This shows that $ f(\tau) =  \c_{\Phi}(\alpha)$, for  
   $\alpha:= h_\topo/ \tau$. Finally, define $\beta_{\Phi}(\alpha):= 
   \gamma_{\Phi}(t_{\tau})$ and note that 
   $\mathcal{P}(\gamma_{\Phi}(t_{\tau})\, I_{\Phi}+ 
 t_{\tau} \varphi)=0$ immediately implies that 
 $P(\beta_{\Phi}(\alpha))= t_{\tau}  \, h_\topo$. 
   Using this and rewriting the above in terms 
  of $\alpha$, the result follows.
 \end{proof}
 
 \section{Multifractal rigidity for NAS}

 For the proof of the main result of this section (Theorem  \ref{NAS}) 
 we require the following proposition. Note that for $u = 
 \delta_{\Phi}$ this result has  been obtained  by 
 Mauldin and Urba\'nski (\cite{MU} Theorem 
 6.1.3). Since it is straight forward to adapt the arguments in \cite{MU}
 to our 
 multifractal situation here,  we will only  give an outline  
 of the proof emphasising the major changes which have to be made.

 \begin{proposition}\label{propext}
	Let $\Phi \in {\mathcal T}_{N}(\Sigma_{d})$ and $u 
	\in \R \setminus \{0\}$ be given, and let  $m_{\Phi,u}$ refer 
	to the    $(u I_{\Phi}-P_{\Phi}(u))$--conformal measure in the 
        measure class 
	of $\mu_{\Phi,u}$.
	Then there exists an open connected  set $W
	\supset X$ such that
	$d\mu_{\Phi,u} /d m_{\Phi,u}$ has a 
	positive real-analytic extension to  $W$. 
\end{proposition}
\begin{proof}(Sketch)
The first step consists of applying Arzel\`a--Ascoli to obtain that
\[ F:C(X)\to C(X), F(g):=\e^{-P_\Phi(u)}\sum_{i\in I}|\phi_i'|^u 
g\circ\phi_i\]
is an almost periodic operator, that is 
$\{F^n(g):n\in \N\}$ is relative compact with respect to the
sup-norm for every $g\in C(X)$ (see  \cite{MU} Lemma 6.1.1).
Also, the Gibbs--property of 
$\mu_{\Phi,u}$ immediately implies that
$F^n(\1)$  is uniformly bounded away from zero and infinity, for 
each $n\in \N$. \\
The second step is to use the above results to show  that there 
exists a unique positive  continuous 
function $\rho:X\to \R^{+}$ such that  (see \cite{MU} Theorem 
6.1.2)
\[F(\rho)=\rho,\; \int \rho \;\dd m_{\Phi,u}=1, \;\mbox{ and } 
\rho|_{\Lambda(\Phi)}=\frac{d \mu_{\Phi,u}}{d m_{\Phi,u}}.\]
The final step is to consider the sequence of functions 
$\left(b_{n}\right)_{n 
\in \N }$, given by  
\[b_n(z):=\sum_{|\omega|=n}|\phi_i '(z)|^u \e^{-nP_\Phi (u)}.\]
One verifies that each $b_{n}$ is  
defined locally on a sufficiently large neighbourhood of each $w \in X$, 
where it is
analytic, uniformly bounded and equicontinuous (see \cite{MU}, proof 
of Theorem 6.1.3).  It  then 
follows  that $(b_{n})$ has a subsequence 
converging to an analytic function which locally extends $\rho$. Since $X$
is compact and simply connected, this provides us with a globally 
defined analytic extention of $\rho$, which is uniformly bounded from above 
and below.
\end{proof}
The following theorem gives the main results of this section. Here, 
the main outcome is that if we have equality up to 
permutation of two Gibbs measures
associated with two points in the Lyapunov spectra of two NAS, then 
the two systems are already bi-Lipschitz equivalent. 
Therefore,  the theorem represents a refinement of the 
Hanus--Urba{\'n}ski rigidity theorem mentioned in the introduction (see also
Corollary \ref{CorNASrigid}).

 \begin{theorem}[Multifractal rigidity  for NAS]\label{NAS}   $ \, $ \\
     Let $\Phi, \Psi \in {\mathcal T}_{N}(\Sigma_{d})$ and 
     $u , v \in  \R \setminus \{0\}$ be given. Then the following 
    three statements are 
 equivalent.
 \begin{itemize}
     \item[  (i)] $\mu_{\Phi,u} \cong \mu_{\Psi,v}$;
     \item[ (ii)] 
     $\Phi \sim \Psi $ and $u=v$;
     \item[(iii)] $I_\Phi\simeq I_\Psi$ and $u=v$.
 \end{itemize}
 Also, the following two 
     statements are 
 equivalent.
 \begin{itemize}
     \item [(iv)] $P_\Phi = P_\Psi$;
     \item[ (v)] $\chi_{\Phi} = \chi_{\Psi}.$
     \end{itemize}
     Furthermore, each of the statements in  {\rm (i) - (iii)} implies the 
     statements in {\rm (iv)} and {\rm (v)}. 
 \end{theorem}
     \begin{proof}
 The implications ``(ii)$\implies$(i)'' ,  ``(iii) $\implies$(i)'' 
  and ``(iii)$\implies$(iv)'', as well as the equivalence of 
  (iv) and (v)  follow exactly   as in the case $\Phi 
  \in {\mathcal T}_{E}(\Sigma_{d})$, and  for this we refer to 
  Theorem \ref{linearrigid} in  Section \ref{Sec:linearrigid}. 
  
  On the basis of the assumption that ``(i)$\implies$(ii)''  holds,
   the implication  ``(i)$\implies$(iii)'' can be obtained as 
   follows.  Assume that $\mu_{\Phi,u} \cong \mu_{\Psi,v}$.  
   We then have 
  $u I_\Phi\simeq v I_\Psi + c$, for some constant $c$. Also, since 
  ``(i)
   $\implies$ (ii)''  holds, we have that $u=v$ and $\Phi \sim \Psi$.
  It hence follows that 
  $I_\Phi - I_\Psi \simeq  c$
  and $\delta_{\Phi}=\delta_{\Psi}$. Consequently, $0= 
  P_{\Phi}(\delta_{\Phi}) - P_{\Psi}(\delta_{\Psi})  = c 
  \delta_{\Phi}$.
  Since $\delta_{\Phi} \neq 0$, this implies that $c=0$, and hence the statement in
   (iii) follows.
  
  It remains to show that ``(i) $\implies$ (ii)''. For this note 
  that by applying a suitable permutation if necessary, we can
  assume without loss of generality that $\mu_{\Phi,u} = \mu_{\Psi,v}$.
  Let $h:\Lambda(\Phi) \to \Lambda(\Psi)$ refer to  
the associated measurable boundary 
  correspondence, that is $\psi_{i}\circ h = h \circ \phi_{i}$ for all 
  $i \in I$. The aim is to show that there exists an open 
  neighbourhood  of $\Lambda(\phi)$ such that 
  $h$ extends to a real-analytic  map on this neighbourhood. For ease of notation
  we will not distinguish between  measures on $\Sigma_{d}$ and
  measures on the corresponding limit sets arising from  representations 
  of $\Sigma_{d}$. For $\omega \in I^{n}$, 
  define
  $J_{\phi_\omega, u}:= d \mu_{\Phi,u} \circ \phi_{\omega} 
  /d \mu_{\Phi,u}$,
  and similar  $J_{\psi_\omega, v}$ for the system $\Psi$.
  We then have for each $i \in I$,
  \begin{eqnarray*} 
      J_{\phi_{i}, u}  & = &  \frac{d \mu_{\Phi,u} \circ \phi_{i}}{d \mu_{\Phi,u}}
      = 
      \frac{d \mu_{\Psi,v} \circ h  \circ  \phi_{i}}{d \mu_{\Phi,v}}
     = \frac{d \mu_{\Psi,v} \circ  \psi_{i} \circ h  }{d 
      \mu_{\Psi,v} \circ h }  = 
      \frac{d \mu_{\Psi,v} \circ  \psi_{i}  }{d 
	   \mu_{\Psi,v} }\circ h \\
	 &  = &  J_{\psi_{i}, v} \circ h.
     \end{eqnarray*}
	   On the other hand, we have
\begin{eqnarray*}  J_{\psi_{i},v} & = & 
	      \frac{d \mu_{\Psi,v} \circ \psi_{i}}{d 
	  m_{\Psi,v} \circ \psi_{i}}  \, \, \frac{d m_{\Psi,v} 
	  \circ \psi_{i}}{d 
	  m_{\Psi,v}}  \, \, \frac{d m_{\Psi,v}}{d 
	  \mu_{\Psi,v}}\\
	  &=& \frac{d \mu_{\Psi,v} }{d 
	  m_{\Psi,v} }\circ \psi_{i}  \, \, \frac{d m_{\Psi,v} 
	  \circ \psi_{i}}{d 
	  m_{\Psi,v}}  \, \, \left(\frac{d \mu_{\Psi,v}}{d 
	  m_{\Psi,v}}\right)^{-1}.
\end{eqnarray*}
Also, since  $m_{\Psi,v} $ is the $(v I_{\Psi}-P_{\Psi}(v))$--conformal 
measure in the measure class 
	of $\mu_{\Psi,v}$,
	  we have
 \[ \frac{d m_{\Psi,v} \circ \psi_{i}}{d 
	   m_{\Psi,v}} = \left| \psi'_{i}\right|^{v}  \, \e^{-P_{\Psi}(v)}
	   , \hbox{ for all } i \in I.
\]
 Now, the conformality condition in the definition of a CS immediately 
 gives that $\left| \psi'_{i}\right|^{v}$ 
	     has a real-analytic extension to an open neighourhood 
	     of $X$. Hence, by combining these observations with  Proposition \ref{propext}, 
	     it follows that there exist  $W \supset X $ such that 
             $J_{\psi_{i}, v}$ has a 
	     real-analytic extension  $\tilde{J}_{\psi_{i}, v}$ to $W$. In the 
	     same way we obtain a 
	     real-analytic extension  $\tilde{J}_{\phi_{i}, u}$ for the system 
	     $\Phi$.
	     Next, note that since $\Psi \in {\mathcal T}_{N}(\Sigma_{d})$, 
	     there exists $j \in I$ such that $\tilde{J}_{\psi_{j}, v}$ is not 
	     equal to a constant.
	     Since $\tilde{J}_{\phi_{j}, u}  = \tilde{J}_{\psi_{j}, v}  \circ h$, 
             the same holds for 
	     $\tilde{J}_{\phi_{j}, u}$  (note, $h$ is defined on the perfect set 
	     $\Lambda(\Phi)$). In particular, the set of zeros of 
	     $\tilde{J}_{\phi_{j}, u}'$, and
	     $\tilde{J}_{\psi_{j},v}'$ respectively,  can not have points of accumulation in
	     $X$, and $Y$ respectively. Therefore, there exists $x \in \Lambda(\Phi)$ 
	    such that $\tilde{J}_{\phi_{j}, u}'(x)  \neq 0$ and  
	    $\tilde{J}_{\psi_{j}, v}'(h(x))  \neq 0$.
	    This implies that there exists an inverse branch 
	    $\tilde{J}_{\psi_{j}, v}^{-1}$ 
	   which is analytic  in a neighbourhood of $\tilde{J}_{\phi_{j}, u}(x)$ such that 
	   $\tilde{J}_{\psi_{j}, v}^{-1}\left(\tilde{J}_{\phi_{j},u}(x)\right)=h(x)$. By choosing 
	   a neighbourhood $W' \subset X$ of $x$  sufficiently small, we obtain 
	   that
	   $\tilde{J}_{\psi_{j}, v}^{-1}\circ \tilde{J}_{\phi_{j}, u}$ is well-defined 
	   and bijective 
	   on $W'$, and  
	   \[ \tilde{J}_{\psi_{j}, v}^{-1}\circ \tilde{J}_{\phi_{j}, u} 
	   (y)=h(y),  \hbox{  for 
	   all  }  y \in W' \cap \Lambda(\Phi).\]
	   It now follows that there exists $\omega \in I^{n}$, for some $n \in 
	   \N$, such that $\phi_{\omega}(X) \subset W'$. Hence, 
	   there exists $W''\supset X$ on which
	   $\psi_{\omega}^{-1} \circ 
	   \tilde{J}_{\psi_{j}, v}^{-1}\circ \tilde{J}_{\phi_{j},u} \circ 
	   \phi_{\omega}$ is real-analytic and 
	   such that $\psi_{\omega}^{-1} \circ 
	   \tilde{J}_{\psi_{j}, v}^{-1}\circ \tilde{J}_{\phi_{j},u} \circ 
	   \phi_{\omega}$ coincides with $h$ on 
	   $W'' \cap \Lambda(\Phi)$.
	   \end{proof}
	   
The following corollary is an immediate consequence of the previous 
theorem. We remark that the fact that $\mu_{\Phi} \cong \mu_{\Psi}$ 
implies that the two Lyapunov spectra coincide is  somehow 
characteristic for non-essentially affine systems. Namely, as we will 
see in Section \ref{EAS}, in this 
respect essentially affine systems behave  rather 
different. Also, note that the equivalence of (i) and (ii) is 
precisely the content of the Hanus--Urba{\'n}ski rigidity theorem.
 \begin{corollary}\label{CorNASrigid}
For $\Phi, \Psi \in {\mathcal T}_{N}(\Sigma_{d})$, the following 
statements are equivalent.
\begin{itemize}
 \item[  (i)] $\mu_{\Phi} \cong \mu_{\Psi}$; 
\item[ (ii)] $\Phi \sim \Psi $.
\end{itemize}
In particular, we also have 
\[ \mu_{\Phi} \cong \mu_{\Psi} \implies \c_{\Phi} = \c_{\Psi}.\]
\end{corollary}	
 {\bf Remark:} Recently, it has been shown in \cite{PW}  that
for cocompact Fuchsian groups
the pressure function is {\em not} a complete invariant of isometry, 
that is equality of the 
pressure functions of two isomorphic cocompact Fuchsian groups
does not necessarily imply that the two associated Riemann surfaces 
are isometric.  This result suggests that one might expect that for 
two systems
$\Phi, \Psi \in  {\mathcal T}_{N}(\Sigma_{d})$ we have  that  
$P_{\Phi} = P_{\Psi}$ does not necessarily imply 
$\Phi \sim \Psi$. However, the argument in \cite{PW} relies on 
Buser's constructive example of  isospectral but  non-isometric,  
compact Riemann surfaces (see \cite{Buser}), and it is currently not 
clear (at least to the authors) how to 
adapt this construction to the situation of a NAS.

\section{Multifractal rigidity and flexibility for EAS}

\subsection{Deformation spaces for EAS}\label{EAS}
We require the following elementary facts about how to switch forward 
and backward 
between two given essentially affine iterated function systems. 
\begin{lemma}\label{hoelder}
Let $\Phi =(\phi_i : X \to \hbox{Int} X \, |  \, i \in 
    I), \Psi =(\psi_i : Y \to \hbox{Int} Y \, |  \, i \in 
    I) \in {\mathcal T}_{E}(\Sigma_{d})$ 
be given. Then there exists a H{\"o}lder 
continuous homeomorphism $h:\Lambda(\Phi) \to \Lambda(\Psi)$ such 
that
\[ \phi_i\circ h=h\circ \psi_{i},  \hbox{  for all $i\in I$}.\]
Moreover, if $\phi_i '=\psi_{i} '$  for all $i\in I$, 
then $h$ is bi-Lipschitz. 
\end{lemma}

\begin{proof}
    Let $\Phi, \Psi $ be given as stated in the lemma. Without loss of generality we can assume that $X=Y=[0,1]$ and 
    that both systems are affine. 
For each $n \in \N$, we define a piecewise linear map $h_n$ by induction as 
follows. For $i \in I$ let 
$I_i:=\phi_i(\Lambda (\Phi))$ and 
$J_i:=\psi_i(\Lambda(\Psi))$, 
and define $h_{0,i}:\hbox{Conv}(I_i)\to \hbox{Conv}(J_i)$  to be the uniquely 
determined linear 
surjection from $I_{i}$ onto $J_{i}$, 
where $\hbox{Conv}$  refers to 
the convex hull. The map $h_0  := \sum_{i=1}^{n} 
h_{0,i}$
is piecewise linear and maps $\bigcup_{i\in I} 
\hbox{Conv}(I_i)$ onto $\bigcup_{i\in I} \hbox{Conv} (J_i)$. 
Similarly, for each $\omega 
\in I^{n}, i \in I$ and $n \in \N$, let $h_{\omega,i}$ be the uniquely 
determined linear 
surjection which maps $\hbox{Conv}(\phi_{\omega i}(\Lambda(\Phi)))$
onto $\hbox{Conv}(\psi_{\omega i} (\Lambda(\Psi)))$. 
Hence, $h_{n}:= \sum_{\omega \in I^{n} }
\sum_{i\in I} h_{\omega,i}$ is a piecewise linear surjection mapping 
$\bigcup_{\omega \in I^{n}} \bigcup_{i\in I} 
\hbox{Conv}(\phi_{\omega}(I_i))$ onto $\bigcup_{\omega \in I^{n}}  \bigcup_{i\in I} 
\hbox{Conv} (\psi_{\omega} (J_i))$.
 Also, one readily verifies that $h_n$ converges uniformly to a continuous function 
$h:=\lim_{n\to\infty} h_n$. The fact that $h$ is H{\"o}lder continuous 
with H{\"o}lder exponent $s:=\min 
\left\{\log(\psi_i') / \log(\phi_i'):i\in I\right\}$ can  be seen as  follows. 
Let $x= \pi_{\Phi} (x_{1} x_{2} \ldots)$ and $ y =\pi_{\Phi} 
(y_{1} y_{2} \ldots)$ be two 
distinct elements of $ \Lambda(\Phi)$. If $x_{1}\neq y_{1}$ then  the 
assertion follows immediately, and hence we can assume without loss of 
generality that $x_{1}=y_{1}$. 
Then there exists
a smallest $n \in {\N}$ such that $x_{n+1}\neq y_{n+1}$ and 
$x_{i}=y_{i}$, for all $1\leq i \leq n$. The open set condition gives 
that there exists $c>0$ 
such that $|x-y| \geq c \prod_{i=1}^{n} \phi_{x_{i}}' $. Using this, 
we obtain
    \[|h(x)-h(y)|  \leq  \prod_{i=1}^{n} \psi_{x_{i}}' 
     \leq  \prod_{i=1}^{n} \phi_{x_{i}}'^{s} 
      =  \frac{1}{c^{s}} \left( c \prod_{i=1}^{n} \phi_{x_{i}}' 
     \right)^{s} \leq \frac{1}{c^{s}} |x-y|^{s}.
     \]
     The remainder of the proposition is now straight forward.
\end{proof}
Note that we necessarily have that each equivalence class in ${\mathcal
T}_{E}(\Sigma_{d})/ \sim$ contains an affine fractal representation.  
Also, note that each affine
fractal representation $\Phi=
(\phi_i : X \to \hbox{Int} X \, | \, i
\in I)$ can be parameterised by 
its {\em contraction
rate vector} $(\phi_{1}',\ldots,\phi_{d}')$, and the
previous lemma  shows that this vector has to be unique
up to permutations of its entries. 
Therefore, as an immediate
consequence of the previous lemma we obtain the following.

\begin{proposition}\label{corEAS}
There exists a canonical bijection from ${\mathcal 
T}_{E}(\Sigma_{d})/ \sim$ onto 
\[ \{(\lambda_{1},
\ldots,\lambda_{d} ) \in  ({\R}^{+})^{d}: 
\sum_{i=1}^{d}\lambda_{i} \leq 1\}/\Pi_d.\]  
Here, $\Pi_{d}$ refers to the group of  permutations of the 
elements in $I$.
\end{proposition}

\subsection{Multifractal rigidity for EAS}\label{Sec:linearrigid}
  
The goal of this section is to study rigidity  for essentially affine 
iterated function systems.
We show that  for these systems  one can only obtain a  multifractal version of 
Sullivan's purely measurable 
rigidity theorem which is significantly weaker than  the one for the  
non-essentially 
affine situation which we obtained in
the previous section.

The following theorem states the main result of this section.
In there it is  shown that in 
the EAS setting
there is a 1-1 
correspondence between the space of pressure functions
and the moduli space ${\mathcal T}_{E}(\Sigma_{d})/ \sim$. 
Also, the theorem in particular gives that for essentially 
affine systems 
equivalence  of $\mu_{\Phi}$ and $\mu_{\Psi}$ alone does in 
general  not  imply
that the  pressure functions of the systems coincide. In fact, as we 
will see in Section \ref{Subsec:MultiFlex}, this 
will only be  the case if the two systems are equivalent. 
Clearly, this can be seen as a first instance exhibiting the 
difference between the essentially 
affine  and the non-essentially affine settings.

\begin{theorem}[Multifractal rigidity for EAS]\label{linearrigid}  $ \, $ \\
    For $\Phi, \Psi \in {\mathcal T}_{E}(\Sigma_{d})$ non-degenerate, the following 
    statements are 
equivalent.
\begin{itemize}
    \item[  (i)] $\mu_{\Phi,u} \cong \mu_{\Psi,u}$ and  $
P_{\Phi}(u)=P_{\Psi} (u)$, for some $u \in \R \setminus \{0\}$;
    \item[ (ii)] $\Phi \sim \Psi $;
    \item[(iii)] $I_\Phi\simeq I_\Psi$;
    \item [(iv)] $P_\Phi = P_\Psi$;
    \item[ (v)] $\chi_{\Phi} = \chi_{\Psi}.$
    \end{itemize}
\end{theorem}

\begin{proof}
  Let $\Phi =(\phi_i : X \to \hbox{Int} X \, |  \, i \in 
	I), \Psi =(\psi_i : Y \to \hbox{Int} Y \, |  \, i \in 
	I) \in {\mathcal 
T}_{E}(\Sigma_{d})$ be two given non-degenerate systems. 
	
	``(i)$\implies$(ii)'':
	Suppose that $\mu_{\Phi,u}=\mu_{\Psi,u}$, for some $u \in 
	\R \setminus \{0\}$ . We then have
	 for each $n \in \N$ and $\omega \in
	I^{n}$,
	\begin{eqnarray*}
	|\phi_\omega (\Lambda(\Phi))|&\asymp & (\mu_{\Phi,u}\circ 
        \pi_{\Phi}^{-1}(\phi_\omega 
	(\Lambda(\Phi))))^{1/u}  \e^{nP_{\Phi}(u)/u}\\
	 & =  & (\mu_{\Psi,u}\circ \pi_{\Psi}^{-1}(\psi_\omega 
	(\Lambda(\Psi)))^{1/u}  \e^{nP_{\Psi}(u)/u} \asymp  |\psi_\omega (\Lambda(\Psi))|.
	\end{eqnarray*}
	We can now proceed similar as in Proposition \ref{hoelder} to 
	build up a bi-Lipschitz map $h: \Lambda(\Phi) \to \Lambda(\Psi)$ as the limit of 
	piecewise linear surjections.  (Note that  the 
	existence of $h$  can be obtained  alternatively by applying Theorem 
	2.2 in \cite{HU}).

	``(ii)$\implies$(i)'': Suppose that
$\Phi \sim \Psi$, and note that  a bi-Lipschitz conjugation does 
not alter the pressure 
function. Hence, similar as in the previous case, 
 we obtain for each $n \in \N$, $\omega \in
I^{n}$ and $u \in 
	\R$,
\begin{eqnarray*}
\mu_{\Phi,u}\circ \pi_\Phi^{-1}(\phi_\omega (\Lambda(\Phi)))&\asymp& |\phi_\omega
(\Lambda(\Phi))|^{u} \e^{-nP_{\Phi}(u)}\asymp |h(\phi_\omega 
(\Lambda(\Phi)))|^{u} \e^{-nP_{\Phi}(u)}\\
& \asymp & |\psi_\omega
(\Lambda(\Psi))|^{u}  \e^{-nP_{\Psi}(u)} \asymp  \mu_{\Psi,u} \circ \pi_{\Psi}^{-1}
(\psi_\omega (\Lambda(\Psi))).
\end{eqnarray*}
Therefore, using the ergodicity of $\mu_{\Phi,u}$ and $ \mu_{\Psi,u}$, it
follows that $\mu_{\Phi,u}=\mu_{\Psi,u}$.

``(i)$\iff$(iii)'':
This is an immediate consequence of the fact that  $\mu_{\Phi,u}$ and
$ \mu_{\Psi,u}$ are Gibbs 
measures for the potential $u I_\Phi - P_{\Phi}(u)$, and $u I_\Psi - 
P_{\Psi}(u)$ 
respectively.

``(iii)$\implies$(iv)'':
This follows from the 
definition of the pressure function. 

``(iv)$\iff$(v)'':
This follows since $P_\Phi$ and $\chi_\Phi$ 
are a Legendre transform pair. 

``(iv)$\implies$(iii)'':
Suppose that $P_\Phi = P_\Psi$, and let $\Phi_a$ 
and $\Psi_a$ be the affine fractal representations within 
the equivalence classes $[\Phi], [\Psi] \in {\mathcal 
T}_{E}(\Sigma_{d})/ \sim$. Also,  
let $(\lambda_1,\ldots ,\lambda_d)$ and
$(\rho_1,\ldots ,\rho_d)$ refer to the
contraction rate vectors  associated with $\Phi_{a}$, and $\Psi_{a}$ 
respectively.  Using the fact that (i) implies (iv), we obtain 
\[P_{\Phi_a}=P_\Phi=P_\Psi=P_{\Psi_a}.\]
Since for affine systems the pressure function at $u$ is equal to the 
logarithm of the sum of the contraction rates raised to the power $u$,  
it follows that
\[\log \sum_{i=1}^{d}\lambda_i^u=P_{\Phi_a}(u)=P_{\Psi_a}(u)=\log 
\sum_{i=1}^{d} 
\rho_i^u, \hbox{  for all  } u\in \R.\]
We can now employ a finite inductive argument as follows.
Let the $\lambda_{i}$ and $\rho_{i}$ be ordered by their sizes
such that $\lambda_{i_{1}} \geq \lambda_{i_{2}} \geq \ldots \geq 
\lambda_{i_{d}}$ and $\rho_{j_{1}} \geq \rho_{j_{2}} \geq \ldots \geq 
\rho_{j_{d}}$. Since $\sum_{i=1}^{d}\lambda_i^u=
\sum_{i=1}^{d} 
\rho_i^u $, it follows

\[
\left(\frac{\lambda_{i_{1}}}{\rho_{j_{1}}} 
\right)^u=\frac{1+\sum_{m=2}^{d}(\rho_{j_{m}}/\rho_{j_{1}})^u}{1+
\sum_{m=2}^{d}(\lambda_{i_{m}}/\lambda_{i_{1}})^u}.
\]
Since for each
$u \geq 0$ the right hand side in the latter equality lies between $1/d$ and 
$d$, we deduce, by letting $u$ tend to infinity,  that the assumption 
$\lambda_{i_{1}}\neq\rho_{j_{1}}$  gives rise to an immediate contradiction. Hence, we
have  that $\lambda_{i_{1}} =\rho_{j_{1}}$. For the inductive step 
assume that for some $k \in I$ we have $\lambda_{i_{n}}= \rho_{j_{n}}$, for all $n\in 
\{1,\ldots,k\}$. We then have $\sum_{m=k+1}^{d}\lambda_{i_{m}}^u=
\sum_{m=k+1}^{d} 
\rho_{j_{m}}^u $, and hence
\[
\left(\frac{\lambda_{i_{k+1}}}{\rho_{j_{k+1}}} 
\right)^u=\frac{1+\sum_{m=k+2}^{d}(\rho_{j_{m}}/\rho_{j_{k+1}})^u}{1+
\sum_{m=k+2}^{d}(\lambda_{i_{m}}/\lambda_{i_{k+1}})^u}.
\]
As above, the right hand side in the latter equality lies between $1/d$ and 
$d$, and hence, by letting $u$ tend to infinity,  we get an immediate
contradiction to the assumption 
$\lambda_{i_{k+1}}\neq\rho_{j_{k+1}}$.  This shows 
that the contraction rate vectors $(\lambda_1,\ldots ,\lambda_d)$ and $(\rho_1,\ldots 
,\rho_d)$ coincide up to a permutation. Combining this observation 
with the fact that (i) implies (iii), it follows that
\[I_\Phi \simeq I_{\Phi_a}=I_{\Psi_a}\simeq I_\Psi.\]
This completes the proof of the theorem.

\end{proof}

The following corollary is an immediate consequence of the previous 
theorem. Note that a comparison of  the statement in here with  
Corollary \ref{CorNASrigid} (see also Theorem \ref{NAS}) clearly 
shows in which respect  essentially 
affine systems have to  be considered as being less rigid than non-essentially 
affine systems.  Also, we remark that  it is straight forward to 
incorporate the degenerate cases. 
\begin{corollary}\label{CorEASrigid}
    For $\Phi, \Psi \in {\mathcal T}_{E}(\Sigma_{d})$, the following 
	statements are 
    equivalent.
    \begin{itemize}
	\item[  (i)] $\mu_{\Phi} \cong \mu_{\Psi}$ and  $
    \delta_{\Phi}=\delta_{\Psi}$;
	\item[ (ii)] $\Phi \sim \Psi $.
	\end{itemize}
Moreover, we have 
	\[ \mu_{\Phi} \cong \mu_{\Psi} \hbox{ and } \delta_{\Phi} = 
	\delta_{\Psi} \implies \c_{\Phi} = \c_{\Psi}.\]
    \end{corollary}
\subsection{Multifractal flexibility for EAS and applications to 
Lyapunov spectra}\label{Subsec:MultiFlex}
As shown in Theorem \ref{linearrigid}, if  for two essentially affine 
systems $\Phi$ and $ \Psi$ we have  that $\mu_{\Phi} \cong \mu_{\Psi}$,
then this does not necessarily imply that the two systems are 
equivalent, nor that their  pressure functions 
coincide.  This naturally raises the question of what can be said about
the pressure functions in case  
$\mu_{\Phi} \cong \mu_{\Psi}$ and $\delta_{\Phi} \neq \delta_{\Psi}$. 
The following theorem  gives a complete answer to this question. 
\begin{theorem}[Multifractal flexibility for EAS]
    \label{Thm:Multiflex}  $ \, $ \\
    For $\Phi, \Psi \in {\mathcal T}_{E}(\Sigma_{d})$ and $u,v\in \R\setminus \{0\}$, the following 
   three statements are equivalent.
\begin{itemize}
    \item[  (i)]  $\displaystyle \mu_{\Phi,u} \cong \mu_{\Psi,v}$;
    \item[(ii)] $\displaystyle I_\Phi\simeq \frac{v}{u} I_\Psi + \frac{P_{\Phi}(u)-P_{\Psi}(v)}{u}$;
    \item[(iii)] $\displaystyle P_{\Phi}(s) =  P_{\Psi}\left(s\cdot \frac{v}{u}\right)+ s\cdot 
                  \frac{P_{\Phi}(u)-P_{\Psi}(v)}{u}$, for all $s \in 
                  \R$.
\end{itemize}
Furthermore, each of the statements in (i) - (iii) implies
\begin{itemize}
   \item[(iv)] $\a_{\Phi} (u) \, \chi_{\Phi}(\a_{\Phi} (u)) = \a_{\Psi} (v)\, 
\chi_{\Psi}(\a_{\Psi} (v)) .$
 \end{itemize} 
\end{theorem}
	 
\begin{proof}  The equivalence ``(ii)$\iff$(iii)''
    can be obtained by exactly the same means as the equivalence ``(iii)$\iff$(iv)''
  in Theorem \ref{linearrigid}. 
Hence,  it is sufficient to show that ``(i)$\iff$(ii)''.

   ``(i)$\iff$(ii)'':  By using a permutation of the generators if necessary, we can 
    assume without loss of generality that $\mu_{\Phi,u} = \mu _{\Psi,v}$.  
    It is then a standard
     result for Gibbs measure that this is equivalent to $u I_\Phi \simeq v
     I_\Psi + P_{\Phi}(u)-P_{\Psi}(v)$, giving that all three
     statements are equivalent.

To finish the proof, it remains to show that (i) and (ii) implies (iv).
 For this, we have by Proposition \ref{MFF}, 
     \begin{eqnarray*} 
 \a_{\Phi} (u) \, \chi_{\Phi}(\a_{\Phi} (u))  &  =& 
u \a_{\Phi}(u) + P_{\Phi}(u) 
 =  -\int \left(u I_{\Phi}-P_{\Phi}(u)
     \right)  \dd \mu_{\Phi,u} \\ &  =&   -\int \left(v I_{\Psi}-P_{\Psi}(v)\right) 
     \dd \mu_{\Psi,v} 
   = v \a_{\Psi} (v) + P_{\Psi}(v)  \\ &  =& \a_{\Psi} (v)\, 
\chi_{\Psi}(\a_{\Psi} (v)) . \end{eqnarray*}
\end{proof}
For the special case in which  $u=\d_\Phi$  and $v=\d_\Psi$, the 
previous theorem has the following immediate implication. 
\begin{corollary}  $ \, $ \\
    For $\Phi, \Psi \in {\mathcal T}_{E}(\Sigma_{d})$, the following 
    statements are 
equivalent.
\begin{itemize}
    \item[  (i)]  $\mu_{\Phi} \cong \mu_{\Psi}$;
    \item[(ii)] $I_\Phi\simeq  \d_{\Psi}/\d_{\Phi} \cdot  I_\Psi$;
    \item[(iii)] $P_{\Phi}(s) = P_{\Psi}(\d_{\Psi}/\d_{\Phi}  \cdot 
    s)$, for all $s \in \R$.
    \end{itemize}
\end{corollary}
Our next aim is to show that there exist systems which
are not bi-Lipschitz equivalent but
which nevertheless admit multifractal measures  which coincide up
to permutation of the generators. For this note that
using Theorem \ref{Thm:Multiflex}, we have
\[\frac{u}{v}I_\Phi=I_\psi + \frac{P_\Phi(u)-P_\Psi(v)}{v}<
\frac{P_\Phi(u)-P_\Psi(v)}{v}.\]
By monotonicity of the pressure function,  it therefore follows 
\begin{equation*}\label{admissible}
P_\Phi\left(\frac{u}{v}\right)<\frac{P_\Phi(u)-P_\Psi(v)}{v}.
\end{equation*}
This observation motivates the following notion
of admissibility.
\begin{definition}\label{adm}
    Let  $\Phi \in \mathcal{T}_{E}(\S)$ and 
      $u,v,p \in 
      \R$ such that $v \neq 0$ be given.
      The triple $(u,v,p)$ is called {\em $\Phi$-admissible} if and only if
      \[ P_{\Phi}\left(\frac{u}{v}\right) < \frac{P_{\Phi}(u) -p}{v}  
      .\]

      \end{definition}
      
 \begin{proposition}[Flexibility  of Lyapunov spectra for 
 EAS (I)]\label{flexsection} $\,$ \\
 Let a non-degenerate $\Phi \in {\mathcal T}_{E}(\Sigma_{d})$  be 
 given,  and let $(u,v,p)$ be a $\Phi$-admissible 
triple. Then there exists
$[\Psi] \in {\mathcal T}_{E}(\Sigma_{d})/ \sim$ (which is unique up to permutations of the 
generators of $\Psi$) such that \[
\mu_{\Phi,u} \cong \mu_{\Psi,v} \;\mbox{ and } \;p=P_{\Psi}(v).\]
 \end{proposition}
 \begin{proof}
  Without loss of generality we can assume that $\Phi$ is an AS.
  Let $(\lambda_{1},\ldots,\lambda_{d}) \in (\R^{+})^{d}$ be the 
  contraction rate vector associated with $\Phi$, and let  
  $(u,v,p)$ be a given $\Phi$--admissible triple.  
      Then define
      \[ \rho_{n}:= \left(\e^{p-P_{\Phi}(u) } \, \lambda_{n}^{u} 
      \right)^{1/v} , \, \hbox{  for each  } n \in \{1,\ldots,d\}.\]
    An elementary calculation immediately shows that the 
    $\Phi$--admissibility of  $(u,v,p)$ 
      is equivalent to  
     $\sum_{n=1}^{d} \rho_{n} < 1$. Hence, by Corollary \ref{corEAS}
     there exists an affine fractal representation $\Psi=(\psi_{i}: 
     [0,1] \to (0,1)  \, | \, i \in I) \in {\mathcal T}_{E}(\Sigma_{d})$
     whose contraction rate vector is $(\rho_{1},\ldots,\rho_{d})$.
     Next, observe that
     \begin{eqnarray*}
	 u I_{\Phi}-P_{\Phi}(u) \hspace{-2mm} &-& \hspace{-2mm} 
         \left( v I_{\Psi}-P_{\Psi}(v)\right) 
	 = -p+ P_{\Psi}(v)  = -p +\lim_{k \to \infty}
	 \frac{1}{k} \log \left(\sum_{n=1}^{d} \rho_{n}^{v}
	 \right)^{k}\\
	 &=& -p + \log \left(  \e^{p-P_{\Phi}(u)} \sum_{n=1}^{d} 
	 \lambda_{n}^{u}\right) = -P_{\Phi}(u) + \log \sum_{n=1}^{d} 
	 \lambda_{n}^{u} =0.
	 \end{eqnarray*}	 
	 This shows that the potentials $u I_{\Phi}-P_{\Phi}(u) $ and  
     $v I_{\Psi}-P_{\Psi}(v)$ coincide,  and also that 
     $p=P_{\Psi}(v)$. 
     It follows that the Gibbs measures corresponding to these 
     potentials have to be equal up to permutation, that is 
     $\mu_{\Phi, u} \cong 
\mu_{\Psi,v}$. 
 \end{proof}

 We end this section by giving a brief discussion of the extent of
 flexibility of an EAS. For this it is more convenient to work
 with the {\em  moduli space of $\Sigma_{d}$} 
     \[ {\mathcal M}_{E} (\Sigma_{d}) := 
	   {\mathcal T}_{E}(\Sigma_{d}) / \sim,\]
where without loss of generality we always assume that an equivalence 
class in  ${\mathcal M}_{E} (\Sigma_{d}) $ is represented by the unique
affine system contained in it. Now, first note that  there  
clearly always is a 
 trivial  measure--wise overlap between the Lyapunov spectra of two 
 EAS, namely $\mu_{\Phi,0} \cong \mu_{\Psi,0}$ for all $\Phi, \Psi 
 \in {\mathcal T}_{E}(\Sigma_{d})$. As we have seen above, for EAS 
 there also is the possibility of non-trivial overlaps, and we 
 will now see that  these are generically represented by 
 $2$--dimensional 
 submanifolds of ${\mathcal M}_{E} (\Sigma_{d})$.

\begin{definition}
Two  systems $\Phi,\Psi \in\mathcal{M}_{E}\left(\Sigma_{d}\right)$
are called \emph{Lyapunov--related} if and only if there exist $u,v\in
\mathbb{R} \setminus \{0\}$
such that $\mu_{\Phi,u}\cong\mu_{\Psi,v}$.
\end{definition}

\begin{figure}[h]\psfrag{1}{$1$}\psfrag{-1}{$-1$}\psfrag{0}{$0$}
\subfigure
[The moduli space 
$\mathcal{M}_{E}\left(\Sigma_{2}\right)$.]
{\includegraphics[width=0.48\textwidth]{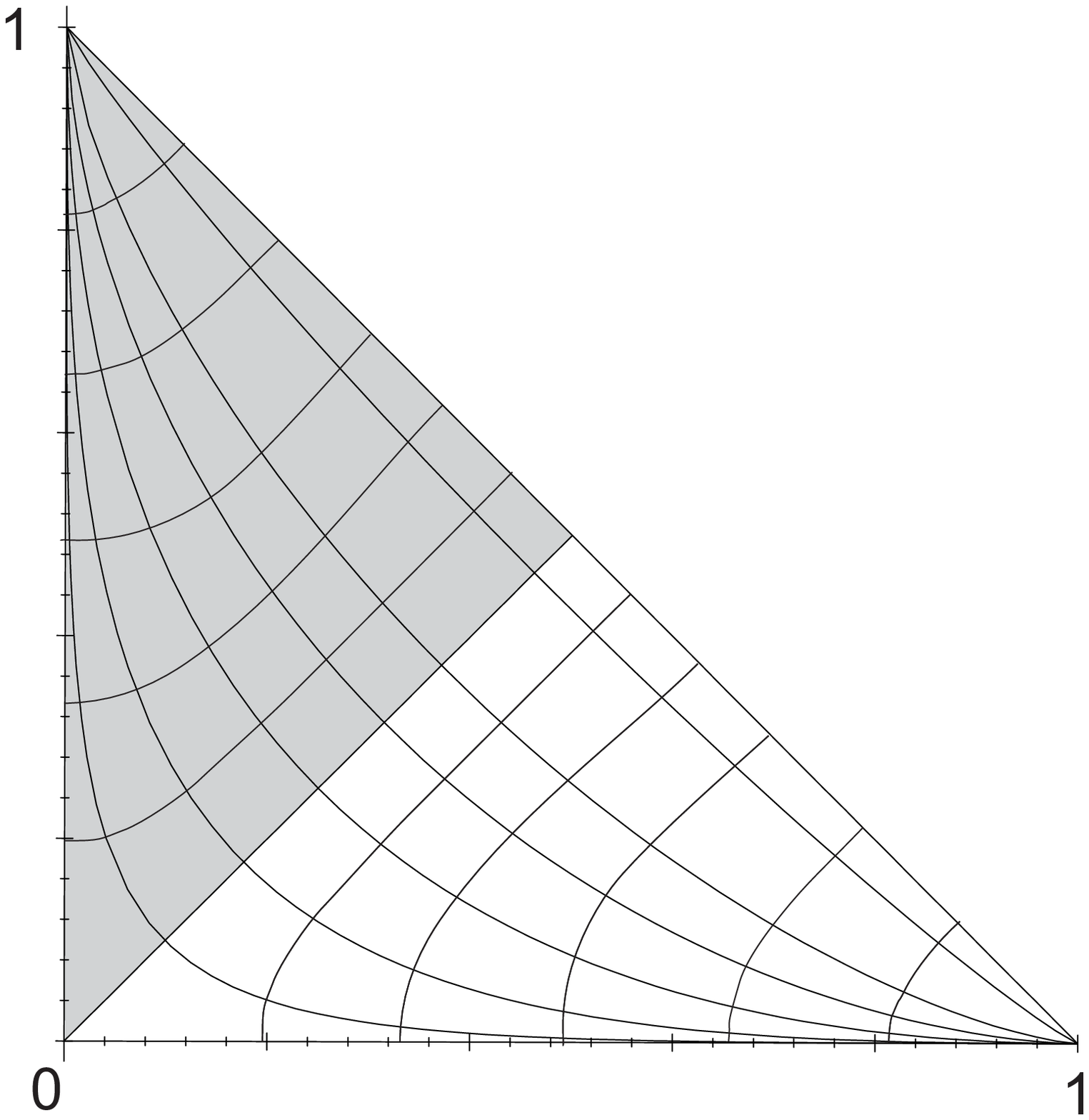}}~~~~
\subfigure
[The harmonized moduli space represented by the disc model.]
{\includegraphics[width=0.48\textwidth]{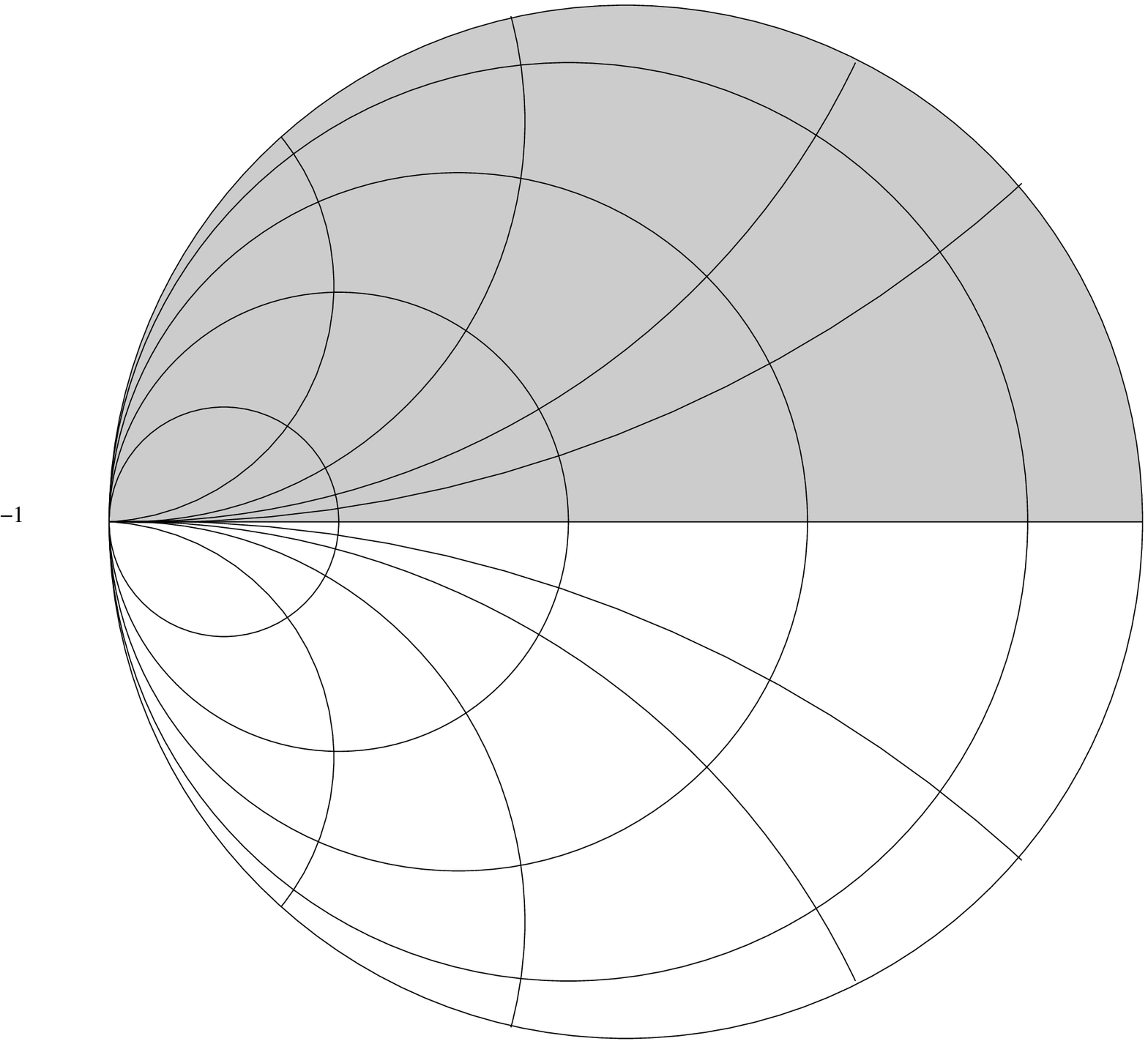}}
\caption{{\it {\bf (a)} The shaded (or alternatively, the non-shaded)
region of the simplex parametrizes the moduli space
$\mathcal{M}_{E}\left(\Sigma_{2}\right):= {\mathcal T}_{E}(\Sigma_{2})
/ \sim$.  The major axes (where at least one generator disappeared)
are not included, whereas the anti-diagonal opposite to the origin
(where the limit set is the whole space $X$) is included.  The
degenerate cases are found on the diagonal.  The lines with endpoints
in (0,1) and (1,0) represent `iso--dimensionals' (i.e. the
Hausdorff-dimension is constant on each of these lines), whereas the
`ortho--dimensionals' (lines orthogonal to the iso--dimensionals) are
the lines of maximal decent of the Hausdorff dimension.  \newline {\bf
(b)} The 'harmonized model' of the moduli space
$\mathcal{M}_{E}\left(\Sigma_{2}\right)$, where the unit intervals on
the major axis are compressed to the singleton $\{-1\} \in {\Sp}^{1}$. 
Here, the iso--dimensionals give rise to the horocyclic
foliation centred at $\{-1\}$, whereas the ortho--dimensionals are
hyperbolic geodesics with one endpoint at $\{-1\}$.}}
\label{fig}
\end{figure}

 Using Theorem \ref{Thm:Multiflex} (ii), we immediately see that 
 if $\Phi$ and  $\Psi$ are  Lyapunov--related, that is 
 $\mu_{\Phi,u}\cong\mu_{\Psi,v}$ for some $u,v\in
\mathbb{R} \setminus \{0\}$, then for each 
$s \in \mathbb{R} \setminus \{0\}$
there exists   $t \in \mathbb{R} \setminus \{0\}$
such that 
$\mu_{\Phi,s}\cong\mu_{\Psi,t}$ (simply choose $t = s \cdot v/u$).
More precisely, we have the following proposition which shows that for a
non-degenerate $\Phi$ the set of systems which are Lyapunov--related
to $\Phi$
forms a $2$-dimensional submanifold of 
$\mathcal{M}_{E}\left(\Sigma_{d}\right)$, whereas if $\Phi$ is 
degenerate 
then this set is a $1$--dimensional submanifold. Note that here, the 
 case $d=2$ appears to be special since it permits only
 exactly two equivalence classes modulo Lyapunov--related,
 namely the diagonal in $\mathcal{T}_{E}\left(\Sigma_{d}\right)$
 and the complement of it in $\mathcal{M}_{E}\left(\Sigma_{d}\right)$
 (see Figure \ref{fig}).
 In all other cases there is a continuum of such equivalence classes.

\begin{proposition}[Flexibility  of Lyapunov spectra for 
 EAS (II)]\label{ending} \hfill
    \begin{itemize}
	\item[ (i)] The `Lyapunov--relation' is an equivalence
relation on $\mathcal{M}_{E}\left(\Sigma_{d}\right)$. 
\item[(ii)] Let $\Phi \in \mathcal{M}_{E}\left(\Sigma_{d}\right)$ be 
given,
and let $(\lambda_{1},\ldots,\lambda_{d})$ be the contraction rate 
vector of $\Phi$. Then the following holds
for the equivalence class $[[\Phi]]$ of $\Phi$ modulo
    the Lyapunov--relation. 
If $\Phi$ is degenerate, then $[[\Phi]] $ is equal to
    \[ \left\{ \Psi \in \mathcal{M}_{E}\left(\Sigma_{d}\right): 
    \rho_{i}
    =  t  , \hbox{ for all } i \in I,  \hbox{ 
    for some } t \in (0,1/d] \right\}.\]
If $\Phi$ is non-degenerate, then $[[\Phi]] $ is equal to
	\[  \left\{ \Psi \in \mathcal{M}_{E}\left(\Sigma_{d}\right): \rho_{i}
	= t \cdot \lambda_{i}^{s}, \hbox{ for all } i \in I,
	\hbox{ 
	    for some } s,t \in \mathbb{R} \setminus
	    \{0\} \right\}. \]
	\end{itemize}
Here, $(\rho_{1},\ldots,\rho_{d})$ refers to the contraction rate 
vector of the system $\Psi$.
     \end{proposition}
     
     \begin{proof}
The assertion in (i) is an immediate consequence of the definition of 
the relation
$\cong$. Furthermore, the first part in (ii) follows since
for degenerate systems the Lyapunov spectrum is trivial. 
For the second part of (ii) we proceed as follows. Let
 $\Phi, \Psi
 \in\mathcal{M}_{E}\left(\Sigma_{d}\right)$
 be two 
non-degenerate systems with contraction rate vector 
$(\lambda_{1},\ldots,\lambda_{d})$, and $(\rho_{1},\ldots,\rho_{d})$
respectively. 
First, if $\Phi$ and $\Psi$ are Lyapunov--related, then
Theorem \ref{Thm:Multiflex}  implies that
there exist $u,v\in\mathbb{R} \setminus \{0\}$ such that \[
\left(\begin{array}{cc}
\log\lambda_{1} & 1\\
\vdots & \vdots\\
\log\lambda_{d} & 1\end{array}\right)\left(\begin{array}{c}
u\\
-P_{\Phi}(u)\end{array}\right)=\left(\begin{array}{cc}
\log\rho_{1} & 1\\
\vdots & \vdots\\
\log\rho_{d} & 1\end{array}\right)\left(\begin{array}{c}
v\\
-P_{\Psi}\left(v\right)\end{array}\right).\]
This implies that  for some $a,b \in \mathbb{R}, a \neq 0$, 
we have\[
\left(\begin{array}{c}
\log\lambda_{1}\\
\vdots\\
\log\lambda_{d}\end{array}\right)=\left(\begin{array}{cc}
\log\rho_{1} & 1\\
\vdots & \vdots\\
\log\rho_{d} & 1\end{array}\right)\left(\begin{array}{c}
a\\
b\end{array}\right) .\]
This settles one direction of the equality. For the reverse direction, assume 
that  
\[
\left(\begin{array}{c}
\log\lambda_{1}\\
\vdots\\
\log\lambda_{d}\end{array}\right)\in\textrm{span}\left(\left(\begin{array}{c}
\log\rho_{1}\\
\vdots\\
\log\rho_{d}\end{array}\right),\left(\begin{array}{c}
1\\
\vdots\\
1\end{array}\right)\right).\]
We  then have that  $I_{\Phi}=vI_{\Psi}+u$, for uniquely determined 
$u,v\in\mathbb{R}$, $v \neq 0$, giving that 
$u= P_{\Phi}(1)-P_{\Psi}(v)$.
Hence, it follows that $I_{\Phi} = vI_{\Psi}+P_{\Phi}(1)-P_{\Psi}(v)$,
which gives $\mu_{\Phi,1}=\mu_{\Psi,v}$.
This shows that  $\Phi$ and $\Psi$
are Lyapunov--related.
\end{proof}

\end{document}